\documentclass[12pt]{amsart}

\usepackage[top=1in, bottom=1in, left=0.9in, right=0.9in]{geometry}
\usepackage{amsmath}
\usepackage{amssymb,amsthm,caption,subcaption}
\usepackage[all]{xy}  
\usepackage{enumerate,graphicx,color}
\usepackage{float}
\usepackage{hyperref}
\usepackage{breakurl}

\theoremstyle{definition}

\theoremstyle{definition}

\title{The software package {\tt S\lowercase{pectral}S\lowercase{equences}}}

\author[Boocher]{Adam Boocher}
\address{Department of Mathematics, University of Utah, Salt Lake City, UT 84112-0090, USA}
\email{boocher@math.utah.edu}

\author[Grieve]{Nathan Grieve}
\address{Department of Mathematics and Statistics, University of New Brunswick, P.O. Box 4400, Fredericton, N.B., Canada, E3B 5A3 }
\email{n.grieve@unb.ca}

\author[Grifo]{Elo\'{i}sa Grifo}
\address{Department of Mathematics, University of Virginia, Charlottesville, VA 22904-4135, USA}
\email{eloisa.grifo@virginia.edu}

\begin{document}
\maketitle

\begin{abstract}
We describe the computer algebra software package {\tt SpectralSequences} for the computer algebra system {\tt Macaulay2}.  This package implements many data types, objects and algorithms which pertain to, among other things,  filtered complexes, spectral sequence pages and maps therein.  We illustrate some of the syntax and capabilities of {\tt SpectralSequences} by way of several examples.
\end{abstract}

\section{Introduction}\label{introduction}
Spectral Sequences play an important role in many areas of mathematics, including Algebraic Topology, Algebraic Geometry and Commutative Algebra.  Here, we describe the software package {\tt SpectralSequences} for the computer algebra system {\tt Macaulay2}.  This package provides tools for effective computation of the pages and differentials in spectral sequences obtained from many kinds of filtered chain complexes.    
We refer the reader to the package's documentation for  precise details regarding the kinds of filtered chain complexes, and the spectral sequences they determine, which we consider.  Further, we mention that  {\tt SpectralSequences version 1.0} is designed to run on {\tt Macaulay2 version 1.9.2} and is available at:
 \url{https://github.com/Macaulay2/  M2/blob/master/M2/Macaulay2/packages/SpectralSequences.m2}.

Spectral sequences, while often very useful, are notoriously cumbersome and difficult to compute.  Indeed, as evidenced by their recursive nature,  making explicit computations by hand is often impossible except when they degenerate quickly.    At the same time, the mere existence of a given spectral sequence can be enough to prove many interesting and important results.  

Before illustrating our package, we mention that we were motivated in part by phrases in the literature, along the lines of:
\begin{itemize}
\item ``There is a spectral sequence for Koszul cohomology which abuts to zero and [provides this formula for syzygies],''  \cite{Green:84};
\item ``Thus the spectral sequence degenerates $\ldots$ [and] $\phi$ is a differential from the $m$th page and other maps are differentials from the first page,'' \cite{Eis:Sch:2009}.
\end{itemize}
We also point out that, in spite of the many classes of spectral sequences which the package {\tt SpectralSequences} can compute, there are still several, including some considered for instance in \cite{Green:84} and \cite{Eis:Sch:2009}, which we would like to compute but which remain out of reach using present techniques and computational power.

In writing this package, our goal was to create a solid foundation and language in {\tt Macaulay2} for working with filtered chain complexes and spectral sequences.  In the class of spectral sequences we can compute, not only can we compute the modules, but we can compute all of the maps as well.  We hope that this package develops further and, as it develops, more examples will be incorporated which will facilitate computations and intuition. Computation and experimentation in algebra in recent decades have 
led to countless conjectures, examples and theorems.  It is our hope that this package will allow for experimentation previously not possible.

As a first easy example, and also to help illustrate one aspect in 
developing our package,   
consider the simplicial complex $\Delta$ on the vertex set $\{x,y,z,w\}$, $x < y < z < w$, with facet description given by $\Delta := \{xyz, wz\}.$
Further, put $F_2 \Delta := \Delta, $ and define simplicial subcomplexes by $F_1 \Delta :=\{xy,w\}$
and $ F_0\Delta := \{x,w\}.$ By considering the reduced chain complexes of the simplicial complexes $F_i\Delta$, over a given  field $k$, we obtain a filtered chain complex
$$
F_\bullet \tilde{\mathcal{C}}_\bullet : 0 = F_{-1} \tilde{\mathcal{C}}_\bullet \subseteq F_0 \tilde{\mathcal{C}}_\bullet \subseteq F_1 \tilde{\mathcal{C}}_\bullet \subseteq F_2 \tilde{\mathcal{C}}_\bullet = \tilde{\mathcal{C}}_{\bullet}.
$$
As it turns out, the spectral sequence
$
E := E(F_\bullet \tilde{\mathcal{C}}_\bullet)
$
determined by this filtered complex 
$F_\bullet \tilde{\mathcal{C}}_\bullet
$
has the property that the map
$$d^2_{2,-1} : E^2_{2,-1} \rightarrow E^2_{0,0}
$$
is an isomorphism of one dimensional $k$-vector spaces.  In fact, it is not difficult to establish this fact directly by hand.  On the other hand, an important first step in developing our  package was for it to successfully compute the map $d^2_{2,-1}$
as well as other similar kinds of examples.  Using the package {\tt SpectralSequences}, we can compute the map $d^2_{2,-1}$
as follows:

\medskip 
\begin{scriptsize}
{\tt 
\begin{verbatim}
i1 : needsPackage "SpectralSequences"; A = QQ[x,y,z,w];

i3 : F2D = simplicialComplex {x*y*z, w*z}; F1D = simplicialComplex {x*y, w}; F0D = simplicialComplex {x,w};

i6 : K = filteredComplex{F2D, F1D, F0D}; E = prune spectralSequence K;

i8 : E^2 .dd_{2,-1}

o8 = | -1 |
              1        1
o8 : Matrix QQ  <--- QQ
\end{verbatim}
}
\end{scriptsize}

In the sections that follow, we briefly describe the structure of our package; we also illustrate it with two examples.  We remark that the package documentation contains many more examples---examples which illustrate how the package can be used to compute, among other things, spectral sequences arising from filtrations of simplicial complexes, triangulations of Hopf fibrations, non-Koszul syzygies and change of rings maps.  

\section{The structure of the package {\tt SpectralSequences} }

The package {\tt SpectralSequences} is able, at least in principle, to compute all aspects of the spectral sequence obtained from a bounded filtered chain complex of finitely generated modules over a finite type $k$-algebra. The actual implementation 
is achieved by first defining a number of auxiliary data structures combined with constructors and other methods to work with these data types.

To use the package {\tt SpectralSequences}, the user must first create a filtered chain complex.  Such filtered chain complexes are represented by the data type {\tt FilteredComplexes}.  The most basic constructor for this type has as input a collection of chain complex maps, whose images determine the given filtered chain complex.  We also provide other methods for creating filtered complexes.  For instance,  the natural filtration induced by truncation can be inputted by the command {\tt filteredComplex}.  For those familiar with the filtered complexes coming from a double complex, we have implemented these in the cases of $\operatorname{Hom}(\mathcal{C}_\bullet,\mathcal C'_\bullet)$ and $\mathcal C_\bullet\otimes \mathcal C'_\bullet$.  Finally, given such a filtered chain complex, represented as an instance of the type {\tt FilteredComplexes}, one uses a constructor associated to the type {\tt SpectralSequence} to create the spectral sequence determined by the given filtered chain complex.

In fact, upon initializing a new spectral sequence, using the type {\tt SpectralSequence}, no calculations are actually performed by the computer.  Rather, calculations are performed using a sort of ``lazy evaluation". There are a number of methods associated to the type {\tt SpectralSequence}.  In brief, for each of the aspects of the spectral sequence, there exists a method which takes, as input,  the spectral sequence represented as an instance of the type {\tt SpectralSequence}.    The output of such methods are either modules, or maps between modules, depending on what is asked by the user.

Finally, the package {\tt SpectralSequences} contains methods which take as input a given spectral sequence represented as type {\tt SpectralSequence} together with a non-negative integer and has as output the modules of the resulting spectral sequence page or the maps between such modules, depending on what is desired by the user.  Such outputs are represented by respective data types {\tt SpectralSequencePage} and {\tt SpectralSequencePageMap}.

\section{Spectral sequences and hypercohomology calculations}

If $\mathcal{F}$ is a coherent sheaf on a smooth complete toric variety $X$, then multigraded commutative algebra can be used to compute the cohomology groups $H^i(X, \mathcal{F})$.  Indeed, if $B \subseteq R$ is the irrelevant ideal of $X$, then the cohomology group $H^i(X, \mathcal{F})$ can be realized as the degree zero piece of the multigraded module $\operatorname{Ext}_R^i(B^{[\ell]}, F)$ for sufficiently large $\ell$; here $B^{[\ell]}$ denotes the $\ell$th Frobenius power of $B$ and $F$ is any multigraded module whose corresponding sheaf on $X$ is $\mathcal{F}$.  
 Given the fan of $X$ and $F$, a sufficiently large power of $\ell$ can be determined effectively. We refer to sections 2 and 3 of \cite{E-M-S} for more details. 
	       
Here we consider the case that $X = \mathbb{P}^1 \times \mathbb{P}^1$ and $\mathcal{F} = \mathcal{O}_C(1,0)$, where $C$ is a general divisor of type $(3,3)$ on $X$.  In this setting, $H^0(C,\mathcal{F})$ and $H^1(C, \mathcal{F})$ are both $2$-dimensional vector spaces.  We can compute these cohomology groups using a spectral sequence associated to a Hom complex. 
	       
We first make the multigraded coordinate ring of 
$\mathbb{P}^1 \times \mathbb{P}^1$, the irrelevant ideal, and a sufficiently high Frobenius power of the irrelevant ideal needed for our calculations.  Also the complex $G$
 below is a resolution of the irrelevant ideal. 

\medskip 

\begin{scriptsize}
\begin{verbatim}
i1 : needsPackage"SpectralSequences";

i2 :  R = ZZ/101[a_0..b_1, Degrees=>{2:{1,0},2:{0,1}}]; -- PP^1 x PP^1 

i3 : B = intersect(ideal(a_0,a_1),ideal(b_0,b_1)); -- irrelevant ideal

o3 : Ideal of R

i4 : B = B_*/(x -> x^5)//ideal; -- Sufficiently high Frobenius power 

o4 : Ideal of R

i5 : G = res image gens B;
\end{verbatim} 
\end{scriptsize}

We next make the ideal, denoted by $I$ below, of a general divisor of type $(3,3)$ on $\mathbb{P}^1 \times \mathbb{P}^1$.  Also the chain complex $F$ below is a resolution of this ideal.

\begin{scriptsize}
\begin{verbatim}
i6 : I = ideal random(R^1, R^{{-3,-3}}); -- ideal of C

o6 : Ideal of R

i7 : F = res comodule I; 
\end{verbatim}
\end{scriptsize}

	       To use hypercohomology to compute the cohomology groups of the 
	       line bundle $\mathcal{O}_C(1,0)$ on $C$ we twist the
	       complex $F$ above by a line of ruling and then 
	       make a filtered complex whose associated spectral
	       sequence abuts to the desired cohomology groups.   This is the complex $K$ below.

\begin{scriptsize}
\begin{verbatim}
i8 : K = Hom(G , filteredComplex (F ** R^{{1,0}}));  

i9 : E = prune spectralSequence K; 

i10 : E^1 

      +-----------------------------------------------+---------------------------------------------+
      | 1                                             | 1                                           |
o10 = |R                                              |R                                            |
      |                                               |                                             |
      |{0, 0}                                         |{1, 0}                                       |
      +-----------------------------------------------+---------------------------------------------+
      |cokernel {-11, 0}  | a_1^5 a_0^5 0     0     | |cokernel {-8, 3} | a_1^5 a_0^5 0     0     | |
      |         {-1, -10} | 0     0     b_1^5 b_0^5 | |         {2, -7} | 0     0     b_1^5 b_0^5 | |
      |                                               |                                             |
      |{0, -1}                                        |{1, -1}                                      |
      +-----------------------------------------------+---------------------------------------------+
      |cokernel {-11, -10} | b_1^5 b_0^5 a_1^5 a_0^5 ||cokernel {-8, -7} | b_1^5 b_0^5 a_1^5 a_0^5 ||
      |                                               |                                             |
      |{0, -2}                                        |{1, -2}                                      |
      +-----------------------------------------------+---------------------------------------------+

o10 : SpectralSequencePage

i11 : E^2; -- output is a mess
\end{verbatim}
\end{scriptsize}

The cohomology groups we want are obtained as follows.

\begin{scriptsize}
\begin{verbatim}

i12 : basis({0,0}, E^2_{0,0}) --  == HH^0 OO_C(1,0)

o12 = {-1, 0} | a_0 a_1 |

o12 : Matrix

i13 : basis({0,0}, E^2_{1,-2}) --  == HH^1 OO_C(1,0)	 

o13 = {-8, -1} | 0               0               |
      {-8, -1} | 0               0               |
      {-8, -1} | 0               0               |
      {-7, -2} | 0               0               |
      {-7, -2} | 0               0               |
      {-7, -2} | 0               0               |
      {-6, -3} | 0               0               |
      {-6, -3} | 0               0               |
      {-6, -3} | 0               0               |
      {-6, -3} | 0               0               |
      {-5, -4} | 0               0               |
      {-5, -4} | 0               0               |
      {-5, -4} | 0               0               |
      {-4, -5} | 0               0               |
      {-4, -5} | 0               0               |
      {-6, -3} | 0               0               |
      {-6, -3} | 0               0               |
      {-6, -3} | 0               0               |
      {-5, -4} | 0               0               |
      {-5, -4} | 0               0               |
      {-5, -4} | 0               0               |
      {-5, -4} | 0               0               |
      {-4, -5} | 0               0               |
      {-4, -5} | 0               0               |
      {-4, -5} | 0               0               |
      {-4, -5} | 0               0               |
      {-4, -5} | 0               0               |
      {-3, -6} | 0               0               |
      {-3, -6} | 0               0               |
      {-3, -6} | 0               0               |
      {-2, -7} | 0               0               |
      {-2, -7} | 0               0               |
      {-2, -7} | a_1^2b_0^4b_1^3 a_1^2b_0^3b_1^4 |

o13 : Matrix
\end{verbatim}
\end{scriptsize}

\section{Seeing cancellations in tensor product complexes}
We now consider an example which is similar to the ones contained in the documentation node {\tt Seeing Cancellations}.  Let $S := k[x_1,\ldots, x_n]$, $I := \langle x_1,\dots, x_n \rangle^2$ and  put $R := S/I$.  By abuse of notation, we write $x_i \in R$ also for the residue class of $x_i \in S$ modulo $I$.   If $K_\bullet$ is the Koszul complex over $R$ on the sequence $(x_1,\ldots, x_n)$ and $F_\bullet$ is a free resolution of $k$ over $R$ then we can form the tensor product complex $D_\bullet = F_\bullet \otimes K_\bullet.$  We now describe these two filtrations and their associated spectral sequences.
\begin{figure}[h!]
\centering
\begin{subfigure}{.5\textwidth}
\xymatrixrowsep{.2in}
\xymatrixcolsep{.2in}
\resizebox{3in}{!}{$$\xymatrix{  
 & \vdots \ar[d]  & \vdots \ar[d] & \vdots \ar[d] & \vdots \ar[d] &  \\
0 & \ar[l] F_0 \otimes K_4 \ar[d]&\ar[l] F_1 \otimes K_4 \ar[d]&\ar[l] F_2 \otimes K_4 \ar[d]&\ar[l] F_3 \otimes K_4 \ar[d]&\ar[l] \cdots \\
0 & \ar[l] F_0 \otimes K_3 \ar[d]&\ar[l] F_1 \otimes K_3 \ar[d]&\ar[l] F_2 \otimes K_3 \ar[d]&\ar[l] F_3 \otimes K_3 \ar[d]&\ar[l] \cdots \\
0 & \ar[l] F_0 \otimes K_2 \ar[d]&\ar[l] F_1 \otimes K_2 \ar[d]&\ar[l] F_2 \otimes K_2 \ar[d]&\ar[l] F_3 \otimes K_2 \ar[d]&\ar[l] \cdots \\
0 & \ar[l] F_0 \otimes K_1 \ar[d]&\ar[l] F_1 \otimes K_1 \ar[d]&\ar[l] F_2 \otimes K_1 \ar[d]&\ar[l] F_3 \otimes K_1 \ar[d]&\ar[l] \cdots \\
0 & \ar[l] F_0 \otimes K_0 \ar[d]&\ar[l] F_1 \otimes K_0 \ar[d]&\ar[l] F_2\otimes K_0 \ar[d]&\ar[l] F_3 \otimes K_0 \ar[d]&\ar[l] \cdots \\
0 & \ar[l]  0 & \ar[l]  0 & \ar[l] 0 & \ar[l] 0 & \ar[l]  & 
}
$$}
\subcaption{A double complex}
\end{subfigure}%
\begin{subfigure}{.5\textwidth}
\centering
\resizebox{3in}{!}{
\xymatrixrowsep{.1in}
$$\xymatrix{  
 & \vdots \ar[d]  & \vdots \ar[d] & \vdots \ar[d] & \vdots \ar[d] &  \\
0 & \ar[l] k(-4)^{n \choose 4} \ar[d]&\ar[l] 0 \ar[d]&\ar[l] 0\ar[d]&\ar[l] 0 \ar[d]&\ar[l] \cdots \\
0 & \ar[l] k(-3)^{n \choose 3}\ar[d]&\ar[l] 0 \ar[d]&\ar[l] 0\ar[d]&\ar[l] 0 \ar[d]&\ar[l] \cdots \\
0 & \ar[l] k(-3)^{n \choose 2} \ar[d]&\ar[l] 0 \ar[d]&\ar[l] 0\ar[d]&\ar[l] 0 \ar[d]&\ar[l] \cdots \\
0 & \ar[l] k(-2)^{n \choose 1}\ar[d]&\ar[l] 0 \ar[d]&\ar[l] 0\ar[d]&\ar[l] 0 \ar[d]&\ar[l] \cdots \\
0 & \ar[l] k(-1)^{n \choose 0} \ar[d]&\ar[l] 0 \ar[d]&\ar[l] 0\ar[d]&\ar[l] 0 \ar[d]&\ar[l] \cdots \\
0 & \ar[l]  0 & \ar[l]  0 & \ar[l] 0 & \ar[l] 0 & \ar[l]  & 
}$$}
\caption{$''E^1$ has degenerated}\label{E1Hor}
\end{subfigure}
\caption{}
\end{figure}
\\ 
\noindent \textbf{The filtrations $'F_\bullet(D)_\bullet$ and $''F_\bullet(D)_\bullet$:} The vertical columns in Figure 1A determine a filtration, whose corresponding spectral sequence $'E = E('F_\bullet(D_\bullet))$ has zeroth page given by 
$
'E^0_{p,q} = F_p \otimes K_q
$
and the maps $'d_{p,q}^0 : F_p \otimes K_q \rightarrow F_p \otimes K_{q-1} 
$ are induced by those of $K_\bullet$.   Analogously, one can consider the filtration $''F_\bullet(D_\bullet)$ determined by the horizontal maps and its spectral sequence $''E$. 

\noindent \textbf{$''E$ degenerates at the first page.}  Since $F_\bullet$ is a resolution, the spectral sequence associated to the horizontal maps degenerates at the first page (see Figure \ref{E1Hor}).  We have
$$''E^\infty =  {''E}^1_{p,q} = \begin{cases} k(-q)^{ n \choose q } & \mbox{if $p = 0$} \\ 0, & \mbox{otherwise}. \end{cases}$$

The second spectral sequence results in an $E^1$ page which will be the Koszul homology of $R$ tensored with $F_\bullet$.  The modules are given by: 
$'E^1_{p,q} = F_p \otimes H^K_q(R)$ where 
$H^K_q(R) \cong \mathrm{Tor}^S_q(R,k) = \bigoplus_j k\left(-j\right)^{\beta_{ij}(R)}.$  Since the $E^\infty$ page consists only of $k(-q)$'s, eventually the $q$th diagonal of $'E^i$ must reduce to only terms generated in degree $q$.  

We can \emph{see} this in our package {\tt SpectralSequences} in any particular example.  For instance, if $S = k[x,y]$ and $I = (x^2,xy,y^2)$, then $'E^1_{p,q}$ is the following complex of $k$-vector spaces.

\begin{figure}[H]
\begin{subfigure}{\textwidth}
\centering
\resizebox{5in}{!}{$$\xymatrix{  
 k \otimes k^2(-3) &\ar[l] k^2(-1) \otimes k^2(-3) &\ar[l] k^4(-2) \otimes k^2(-3)  &\ar[l] k^8(-3) \otimes k^2(-3)  & \ar[l] \cdots\\ 
 k \otimes k^3(-2) &\ar[l] k^2(-1) \otimes k^3(-2) &\ar[l] k^4(-2) \otimes k^3(-2)  &\ar[l] k^8(-3) \otimes k^3(-2)   & \ar[l] \cdots\\ 
 k \otimes k            &\ar[l] k^2(-1) \otimes k           &\ar[l] k^4(-2) \otimes k          &\ar[l] k^8(-3) \otimes k & \ar[l] \cdots
}$$}\caption{$'E^1$}
\end{subfigure}
\vskip .5in
\centering
\begin{subfigure}{.5\textwidth}
\centering
\resizebox{3in}{!}{
$$\xymatrix{  
 k^2(-3) &   k^4(-4) &             k^8(-5)  &   k^{16}(-6)  &  \cdots\\ 
{k^3(-2)} &  k^6(-3) &\ar[ull]   k^{12}(-4)  &\ar[ull]   k^{24}(-5)   & \ar[ull] \cdots\\ 
 \textcolor{blue}{k}            &  \textcolor{blue}{k^2(-1)} &\ar[ull]   \textcolor{blue}{k(-2)}\oplus k^3(-2)  &\ar@[red][ull]   k^8(-3) & \ar[ull] \cdots
}
$$
}
\subcaption{$'E^2$}
\end{subfigure}%
\begin{subfigure}{.5\textwidth}
\centering
\resizebox{3in}{!}{

$$\xymatrix{  
 k^2(-3) &   k^4(-4) &             k^8(-5)  &   k^{16}(-6)  &  \cdots\\ 
0 &  0 &  0  & 0   &  \cdots\\ 
 \textcolor{blue}{k}            &  \textcolor{blue}{k^2(-1)} &   \textcolor{blue}{k(-2)}  &\ar@[red][uulll]   k^2(-3) & \ar[uulll] \cdots
}.
$$
}
\subcaption{$'E^3$}
\end{subfigure}
\end{figure}
We have separated out the $E^\infty$ terms in blue.   Notice that all the remaining (black) terms must disappear at some later page.  In this case, the maps from the $'E_{3,0}$ position are both surjective, which we illustrate below:

\medskip 
\begin{scriptsize}
{\tt 
\begin{verbatim}
i9 : S = ZZ/101[x,y];

i10 : I = ideal(x^2,x*y,y^2);   R = S/I;

i11 : kR = coker vars R;   kS = coker vars S;

i12 : K = (res kS)**R;  F = res(kR,LengthLimit=>6);

i13 : E = prune spectralSequence (K ** filteredComplex F);

i14 : length image ((E_2).dd_{3,0})

o14 : 6

i15 : length image (E_3).dd_{3,0}

o15 : 2
\end{verbatim}
}
\end{scriptsize}

\section*{Acknowledgments}
\vspace{.4cm}
Many people have worked on this package.  The authors especially thank David Berlekamp, Greg Smith, and Thanh Vu for their contributions.  The authors also thank Claudiu Raicu for help with his {\tt PushForward} package and also Dan Grayson for answering questions.  During the final portions of this work, the second author was financially supported by an AARMS postdoctoral fellowship. Finally, the first and second authors respectively thank David Eisenbud and Mike Stillman for suggesting the project to them.

\providecommand{\bysame}{\leavevmode\hbox to3em{\hrulefill}\thinspace}
\providecommand{\MR}{\relax\ifhmode\unskip\space\fi MR }
\providecommand{\MRhref}[2]{%
  \href{http://www.ams.org/mathscinet-getitem?mr=#1}{#2}
}
\providecommand{\href}[2]{#2}

\end{document}